\documentclass[a4paper,12pt,reqno]{amsart}
\usepackage{amssymb}
\usepackage{amsmath}

\def\i{\,\lrcorner\,}
\def\a{\alpha}
\def\b{\beta}
\def\g{\gamma}

\def\vs{\vskip .6cm}

\def\la{\langle}
\def\ra{\rangle}
\def\.{\cdot}
\def\O{\Omega}
\def\n{\nabla}

\def\t{\tilde}
\def\beq{\begin{equation}}
\def\eeq{\end{equation}}
\def\bea{\begin{eqnarray*}}
\def\eea{\end{eqnarray*}}
\def\ba{\begin{array}}
\def\ea{\end{array}}

\def\o{\omega}

\def\L{\Lambda}

\def\r{\end{proof}}

\def\rp{R^+}
\def\I{{\mathcal I}}

\def \RM{\mathbb{R}}

\def \HM{\mathbb{H}}
\def\End{{\rm End}}

\def\d{{\delta}}

\def\raw{\rightarrow}

\def\es{\,\lrcorner\,}

\def\Ric{\mathrm{Ric}}

\def\be{\begin{equation}}

\def\ee{\end{equation}}

\newtheorem{ede}{Definition}[section]

\newtheorem{epr}[ede]{Proposition}

\newtheorem{ath}[ede]{Theorem}

\newtheorem{elem}[ede]{Lemma}

\newtheorem{ecor}[ede]{Corollary}
 
\title{Killing Forms on Symmetric Spaces}
\author{Florin Belgun, Andrei Moroianu and Uwe Semmelmann}
\thanks{The authors are  members of the {\sl European Differential Geometry 
Endeavour} (EDGE), Research Training Network HPRN-CT-2000-00101, 
supported by The European Human Potential Programme;
the first and third authors would like to thank the Centre de
Math{\'e}matiques de l'Ecole Polytechnique for hospitality during the
preparation of this work.}

\address{Florin Belgun\\ Institut f{\"u}r Mathematik \\
 Universit{\"a}t Leipzig\\
Augustusplatz 10-11\\ D-04109 Leipzig, Germany}
\email{belgun@mathematik.uni-leipzig.de}

\address{Andrei Moroianu \\ CMAT\\ {\'E}cole Polytechnique \\ UMR 7640 du CNRS
\\ 91128 Palaiseau \\ France}
\email{am@math.polytechnique.fr}

\address{Uwe Semmelmann\\ Fachbereich Mathematik, Universit{\"a}t
Hamburg\\ Bundesstr. 55,  D-20146 Hamburg, Germany}
\email{Uwe.Semmelmann@math.uni-hamburg.de}

\begin{document}

\begin{abstract}
Killing forms on Riemannian manifolds are differential forms whose covariant
derivative is totally skew--symmetric. We show that a compact
simply connected symmetric space carries a non--parallel Killing $p$--form
($p\ge2$) if and only if it isometric to a Riemannian product
$S^k\times N$, where $S^k$ is a round sphere and $k>p$.
\vs


\end{abstract}

\maketitle

\section{Introduction}

There are two equivalent definitions of Killing vector fields on
Riemannian manifolds. A vector field $X$ is Killing if its local flow
consists of isometries. Equivalently, $X$ is Killing if the covariant
derivative $\nabla X^{\flat}$ of the dual $1$--form $X^{\flat}$ is
skew--symmetric. 

This second definition can be generalized to forms of higher degree as
follows. A $p$--form $u$ is called {\em Killing} if its covariant
derivative is totally skew--symmetric, {\em i.e.} if it exists some
$p+1$--form $\tau$ such that 
$$\n u=\tau.$$
It is easy to check that in that case $\tau$ is necessarily equal to 
$\frac{1}{p+1}d u$. In contrast to Killing $1$--forms, which are just
dual to infinitesimal isometries, there is no geometrical interpretation
of Killing $p$--forms for $p\ge2$. 

The aim of this paper is to show the following

\begin{ath}\label{main}
If a symmetric space $M$ of compact
type carries a non--parallel Killing $p$--form and $p\ge 2$, then the universal
cover $\tilde M$ is either a round sphere, or has a factor isometric
to a round sphere in its de Rham decomposition.
\end{ath}

Together with the fact that a Killing form on a product splits as a product 
of Killing forms on the factors \cite{prod}, we get therefore a complete
description of all Killing forms on locally symmetric spaces of compact type.

This result can be thought of as a generalization of the following
weaker assertion: 

\begin{epr} \label{kspi} A symmetric space admitting real Killing spinors is
locally conformally flat. 
\end{epr}

To see that Theorem \ref{main} implies Proposition \ref{kspi}, we
recall the fact that a manifold $M$ which carry Killing spinors is locally
irreducible (cf. \cite{bfgk}) and using the squaring construction one
can construct non--parallel Killing $p$--forms for some $p\ge 2$ 
starting from a Killing spinor. Thus the universal cover of $M$ has to
be the sphere, so $M$ is conformally flat. 

Of course, one can give a more direct proof of Proposition
\ref{kspi}. If $\Psi$ is a Killing spinor, an immediate calculation
shows that $M$ is Einstein and $W(X,Y)\.\Psi=0$, where $W$ is the Weyl 
tensor, for all $X,Y\in TM$
(cf. \cite{bfgk}). 
Differentiating this relation several times and using the fact that
the Weyl tensor is parallel, we obtain that for every vectors $X$ and
$Y$, the Clifford product of the $2$--form $W(X,Y)$ with every spinor
vanishes, so finally $W=0$. 

The strategy for the proof of Theorem \ref{main} is somewhat similar,
although much more involved. We first show that if $M$ is a Riemannian
product, then at least one of the factors carries a Killing $p$--form, too.
We then interpret Killing $p$--forms as
parallel sections of $\L^pM\oplus\L^{p+1}M$ with respect to some
modified natural connection $\bar\n$ acting on this bundle. The
curvature $\bar R$ of
this connection can be computed explicitly in terms of the Riemannian
curvature of $M$, and the sections of $\L^pM\oplus\L^{p+1}M$ which lie
in the kernel of $\bar R$ define some $\n$--parallel sub--bundle
$E_0\oplus F_0$ of 
$\L^pM\oplus\L^{p+1}M$. The point is that this sub--bundle is not necessarily
$\bar \n$--invariant. By an inductive procedure, one can construct a
sequence of $\n$--parallel sub--bundles $E_{i+1}\oplus F_{i+1}\subset
E_i\oplus F_i$ such that 
$\bar\n (E_{i+1})\subset F_i$ and $\bar\n (F_{i+1})\subset E_i$. 
This sequence is of course stationary,
and defines some $\n$-- and $\bar\n$--parallel sub--bundle $E\oplus
F=E_k\oplus F_k$ of
$\L^pM\oplus\L^{p+1}M$ for some $k$ large enough. 

A tricky argument (which is the core of the paper and is described in
detail in Section \ref{core}) allows one to show that the projection $F$
of $E\oplus F$ onto $\L^{p+1} M$ is either zero, or the whole
space. The first case just says that every Killing form has to be
parallel, while in the second situation it is easy to show that the
Weyl curvature of $M$ has to vanish.

\section{Preliminaries}

Throughout this paper we use Einstein's summation conventions on
double subscripts. Vectors and $1$--forms are identified via the metric. In the
sequel, $\{e_i\}$ will denote a local orthonormal basis of the tangent
bundle, parallel at some point. 

\begin{ede}
A $p$--form $u$ is called a {\it Killing $p$--form} if and only if 
\begin{equation}\label{killing}
\nabla_Xu=
\frac{1}{p+1}X\es du,
\end{equation}
for all vector fields $X$. 
\end{ede}

Let $u$ be a Killing $p$--form. Obviously $X\es \n_X u=0$ for all
vectors $X$ so in particular $\d u=0$.
Let us take the covariant derivative in (\ref{killing}) with respect
to some vector field 
$Y$, wedge with $X$ and sum over an orthonormal basis $X=e_i$ ($X$,
$Y$ and $e_i$ are supposed to be parallel at a point):
\bea e_i\wedge R_{e_i,Y}u&=& e_i\wedge \n_{e_i}\n_Y u-e_i\wedge \n_Y\n_{e_i}u=
e_i\wedge \n_{e_i}(\frac{1}{p+1}Y\es du)-\n_Y d u\\
&=&-\frac{1}{p+1}Y\es e_i\wedge \n_{e_i}du+\frac{1}{p+1}\n_Y du-\n_Y d
u\\
&=&-\frac{1}{p+1}Y\es d(du)-\frac{p}{p+1}\n_Y du=-\frac{p}{p+1}\n_Y
du.
\eea
Thus, denoting by $R^+$ the operator 
$$R^+:TM\to \End(\L^pM,\L^{p+1}M),\qquad R^+(X)u:=e_i\wedge
R_{X,e_i}u,$$
the above equation reads
\begin{equation}\label{killing2}
\nabla_Xdu=
\frac{p+1}{p}R^+(X)u.
\end{equation}
 
Consider the connection $\t\n$ on $\L^pM\oplus\L^{p+1}M$ given by
$$\t\n_X(u,v):=(\n_Xu-\frac{1}{p+1}X\es v,
\n_Xv-\frac{p}{p+1}R^+(X)u).$$
We have just shown that $(u,du)$ is a $\t\n$--parallel section of
$\L^pM\oplus\L^{p+1}M$ for every Killing $p$--form $u$. We will not
discuss here the consequences of this important fact, but rather
refer to \cite{uwe} for details.

Taking the covariant derivative with respect to some vector field $Y$
in (\ref{killing}), skew--symmetrising in $X$ and $Y$ and using
(\ref{killing2}) yields 
\beq\label{k1}R_{X,Y}u=-\frac1p (X\i R^+(Y)u-Y\i
R^+(X)u)\quad\forall X,Y\in TM,
\eeq
which we rewrite as
\beq\label{c0}
(\I\wedge \rp )u=-p\,R\,u,\ \forall u\in\L^pM \ \mbox{Killing form}.
\eeq
Here $\I$ stands for the interior product, $\I(X)u:=X\i u$, and $\rp $ for 
the operator 
defined before, both of which are viewed as a $1$--form with values 
in $End(\L^*M)$. Therefore their exterior product is a $2$--form with
values in  $End(\L^*M)$, and so is the curvature operator $R$. Note that 
$\I$ decreases and $\rp $ increases the degree of the form by $1$.

If $M$ is locally symmetric, then $R$ is parallel, and so is $\rp $ 
($\I$ is always parallel). Thus, taking the
covariant derivative with respect to $Y$ 
in (\ref{killing2}), skew--symmetrising in $X$ and $Y$ and using
(\ref{killing}) yields 
\beq\label{k2}R_{X,Y}du=-\frac1p (R^+(X)Y\es du-
R^+(Y)X\es du)\quad\forall X,Y\in TM,
\eeq
or, equivalently,
\beq\label{conc1}
(\rp\wedge \I)du=-p\,R\,du,\ \forall u\in\L^pM \ \mbox{Killing form}.
\eeq

We end up this section by deriving some useful algebraic relations satisfied
by $R^+$.

\begin{elem} \label{l1} On any manifold we have

\beq \rp\wedge \I=\I\wedge\rp +R,\mbox{ or, equivalently, }2[\rp\wedge \I]=R.
\eeq

\end{elem}
Here we define, for two $1$--forms $A,B$ with values in some algebra bundle
$\O$, their {\it commutator} $[A,B]:TM\otimes TM\raw \O$ by 
$[A,B](X\otimes Y):=A(X)B(Y)-B(X)A(Y)$, and $[A\wedge B]:\L^2(TM)\raw\O$, 
resp. $[A\odot B]:S^2(TM)\raw\O$ denote its skew--symmetric, resp. 
symmetric, part.
\proof
\bea R^+(X)\I(Y) u &=& e_i\wedge R_{X,e_i}(Y\i u)=e_i\wedge
R_{X,e_i}Y\i u +e_i\wedge Y\i R_{X,e_i} u\\
&=&e_i\wedge R_{X,e_i}Y\i u - \I(Y)\i R^+(X)(u)+R_{X,Y}u,
\eea
which, after skew--symmetrization and using the Bianchi identity on
the first term of the right hand side, yields the desired result.
\qed

\begin{ecor} Let $E$ be a parallel (i.e., $\n$--stable) sub--bundle 
of $\Lambda^pM$.
The induced operator $\widetilde{[\rp ,\I]}:TM\otimes TM\raw Hom(E,\L^p M/E)$ is a 
symmetric 2--tensor.
\end{ecor}
\begin{elem}\label{lema2} Let $E$ be as before. Then the induced 3--tensor
$$\widetilde{\I^2\rp }:TM^{\otimes 3}\raw \L^{p-1}M/(\I(E)+\I\rp \I(E))$$
vanishes identically.
\end{elem}
\proof
First note that
$$\I^2\rp =\I\rp \I+\I[\I,\rp ]=\I\rp \I-\frac12\I R+\I[\I\odot\rp ],$$
therefore 
$$\widetilde{\I^2\rp} \equiv \I[\I\odot\rp ] \mbox{ (mod) }\I(E)+\I\rp \I(E).$$
The tensor $\widetilde{\I^2\rp }$ is clearly skew--symmetric in the first two 
arguments, because two interior products anti--commute. On the other hand, 
the previous equation shows that $\widetilde{\I^2\rp }$ is equal to the
co--restriction to $Hom(E,\L^p M/(\I(E)+\I\rp \I(E))$ of the tensor 
$\I[\I\odot\rp ]$, which is,  
by its very definition, symmetric in the last two arguments.

But a $3$--tensor which is symmetric in the last two arguments and 
skew--symmetric in the first two arguments is necessarily zero.
\qed

\section{Form bundles on symmetric spaces}\label{core}

The results in this section could have been stated in terms of
abstract representation theory, but we prefer the more geometric
presentation below.

Let $M$ be a symmetric space with curvature tensor $R$.
We define the following vector bundles on $M$:
$$E_0:=\{u\in\L^pM\ |\ Ru=-\frac1p \I\wedge\rp u\},$$
$$F_0:=\{v\in\L^{p+1}M\ |\ Rv=-\frac1p \rp \wedge \I v\}.$$
We then define inductively the vector bundles
$$E_k:=\{u\in E_{k-1}\ |\ R^+(X)u\in F_{k-1}\quad\forall X\in TM\},$$
$$F_k:=\{v\in F_{k-1}\ |\ X\i v\in E_{k-1}\quad\forall X\in TM\}.$$
Since $R$ is parallel, we see that $E_k$ and $F_k$ are parallel vector
bundles for every $k$. We denote by 
\beq\label{ef}E:=\underset{k}\cap E_k,\qquad F:=\underset{k}\cap F_k.\eeq
By definition, for every sections $u$ and $v$ of $E$ and $F$
respectively we have 
\beq\label{c1}Ru=-\frac1p \I\wedge \rp u
\eeq
\beq\label{c2}Rv=-\frac1p \rp \wedge \I v
\eeq
\beq\label{c31}R^+(X)u\in F\qquad \hbox{and}\qquad 
X\i v\in E\quad\forall X\in TM
\eeq
Notice that from (\ref{c31}) we get
\beq\label{c3} \I\rp(E)\subset E.\eeq

\begin{elem} Let $k$ be an integer $k\ge1$. For every tangent vectors
  $X_1,\cdots, X_k$, $Y_1,\cdots, Y_k$ and for every 
section $u$ of $E$ we have 
\beq \label{cont}X_1\i \ldots\i X_k\i R^+(Y_1)\ldots R^+(Y_k)u\in E.\eeq
\end{elem}

\proof
The statement claims that $\I^k ({\rp })^k$, as a $2k$--tensor with values in 
$Hom(E,\L^pM)$, takes actually values in $End(E)$.

We use induction on $k$. For $k=1$ the result follows from (\ref{c3}).

{\it Step 1. } First we show (also by induction), using the Lemma \ref{lema2}, that
\begin{equation}\label{etapa} 
\I^l {\rp }(V)\subset \I\rp\I^{l-1}(V)+ I^{l-1}(V),
\end{equation}
for any $l\ge 1$ and any $\nabla$--parallel 
form sub--bundle $V$. Indeed, for $l=2$ this is exactly Lemma \ref{lema2}, and 
for $l>2$, we get from the induction step that $\I^l {\rp }(V)\subset
\I^2\rp\I^{l-2}(V)+ \I(\I^{l-2}(V))$, from which the claim follows
using again Lemma \ref{lema2}.

{\it Step 2. } Let us denote now by $E':=({\rp })^{k-1}(E)\subset\L^{p+k-1}$ the image of 
the last factors in the product $\I^k({\rp})^k$. 
Then $E'$ is a parallel bundle. We have to show 
that $\I^k\rp(E')\subset E$. From (\ref{etapa}) we get $ \I^k\rp(E')\subset 
\I\rp\I^{k-1}(E')+\I^{k-1}(E')$, and from the induction step 
$\I^{k-1}(E')\subset E$, so we get in the end $\I^k({\rp})^k(E)\subset
E+\I\rp(E)\subset E$ by (\ref{c3}).



\qed

\begin{ecor}\label{co1}If there exist $n-p$ 
tangent vectors $Y_1,\ldots,Y_{n-p}$
and a $p$--form 
$u\in E$ such that $R^+(Y_1)\ldots R^+(Y_{n-p})u$ is non--zero,
then $E=\L^pM$.
\end{ecor}

\proof
This follows simply because the different contractions of the
volume form with $n-p$ vectors span $\L^p M$.

\qed

We now examine under which circumstances the hypothesis in the
corollary above can fail, that is, what can one say about $E$ if
$R^+(Y_1)\ldots R^+(Y_{n-p})u=0$ for all $Y_1,\ldots,Y_{n-p}\in TM$
and $u\in E$.

\begin{elem}\label{lem1}
Let $V\subset\L^qM$ be some irreducible summand in the decomposition
of $q$--forms under the holonomy representation.

{\em (i)} If $R^+(X)u=0$ for all tangent vectors $X$ and 
$u\in V$ then the holonomy representation on $V$ is trivial.

{\em (ii)} Suppose that the holonomy representation is
irreducible on $TM$ and that $M$ is not K{\"a}hler. 
If there exists some sub--bundle $W$ of
$\L^{q+1}M$ on which the holonomy
representation is trivial and such that $R^+(X)u\in W$ for every
tangent vector $X$ and $u\in V$, then either $q=n-1$, or $R^+(X)u=0 \
\forall\,X$.  

\end{elem}

\proof (i) Taking the interior product with $X$ and making the sum
over an orthonormal basis yields
$$0=e_i\i R^+(e_i)u=-q(R)u.$$
But $q(R)$, being the Casimir operator of the holonomy group on $V$,
is a non--negative constant. Moreover, this constant is zero if and
only if $V$ is the trivial representation. 

(ii) Let $H$ be the holonomy group of $M$. Since $M$ is irreducible,
it has to be Einstein, and we denote its Einstein constant by $r$.
The operator $R^+$ defines an equivariant map $R^+:TM\otimes
V\to W$. We have two possibilities: either $V$ is isomorphic (as
$H$--representation) to $TM$, or not. In the last case, the real Schur
Lemma shows that the map $R^+$ vanishes. We now suppose that $V$ is
isomorphic to $TM$. In particular, $q(R)$ acts on $V$ by
multiplication with the scalar $r$ (remember that $q(R)=\Ric$ on
$TM\simeq\L^1M$). Let $\{v_{\a}\}$, $\a=1,\ldots,$dim$(W)$ be an
orthonormal basis of $W$. Since the curvature acts trivially on $W$ we
can write for every $u\in V$
\bea \la R^+(X)u, v_\a\ra&=&\la e_i\wedge R_{X,e_i}u,v_\a\ra=
-\la u, R_{X,e_i}(e_i\i v_\a)\ra\\
&=&-\la u,(R_{X,e_i}e_i)\i v_\a\ra=-r\la u,X\i v_\a\ra,
\eea
whence
$$R^+(X)u=\la R^+(X)u,v_\a\ra v_\a=-r\la u,X\i v_\a\ra v_\a.$$
Taking the interior product in this formula yields
\beq\label{p} u=\frac1r q(R)u=-\frac1r e_i\i R^+(e_i)u=e_i\i \la u,e_i\i
v_\a\ra v_\a =\la u,e_i\i v_\a\ra e_i\i v_\a.
\eeq
Since $M$ is not K{\"a}hler, it turns out that every $H$--equivariant
bilinear form on $TM$ is a real multiple of the metric. In particular,
we get 
$$\la X\i v_\a,Y\i v_\b\ra = c_{\a\b}\la X,Y\ra.$$
Taking $X=Y=e_i$ and summing yields
\beq nc_{\a\b}=c_{\a\b}\la e_i,e_i\ra=\la e_i\i v_\a,e_i\i v_\b\ra
=\la e_i\wedge e_i\i v_\a,v_\b\ra =(q+1)\d_{\a\b},
\eeq
so finally 
\beq\label{p1}\la X\i v_\a,Y\i v_\b\ra = \frac{q+1}{n}\d_{\a\b}\la X,Y\ra.\eeq
We take the scalar product with $e_j\i v_\b$ in (\ref{p}) and use
(\ref{p1}) to obtain
$$\la u, e_j\i v_\b\ra=\la u,e_i\i v_\a\ra\la e_i\i v_\a,e_j\i
v_\b\ra= \frac{q+1}{n}\la u, e_j\i v_\b\ra.$$
Since $u$ is non--zero by assumption, this shows that $q=n-1$.

\qed

\begin{ecor}\label{cor2}
Let $M$ be a compact irreducible non--K{\"a}hlerian symmetric space. If 
the bundle $E$ defined in (\ref{ef}) 
is a proper sub--bundle of $\L^pM$ then $R^+(X)u=0$ for all
$X\in TM$ and $u\in E$.
\end{ecor}

\proof Let $k$ be the smallest positive integer such that 
$R^+(Y_1)\ldots R^+(Y_k)u=0$
for all $ Y_1,\ldots Y_k\in TM$ and $u\in E$. Clearly $1\le k\le n-p+1$.
If $E$ is strictly included in $\L^pM$,
Corollary \ref{co1} shows that $k\le n-p$. Let $W$ be the maximal
sub--bundle of $\L^{p+k-1}M$ on which the holonomy group acts
trivially. From Lemma \ref{lem1} (i)
we see that $R^+(Y_1)\ldots R^+(Y_{k-1})u\in W$ for all $ Y_1,\ldots
Y_{k-1} \in TM$ and $u\in E$. 

Suppose that $k\ge2$. Lemma  \ref{lem1} (ii) shows that either 
$R^+(Y_1)\ldots R^+(Y_{k-1})u=0$ for all $ Y_1,\ldots
Y_{k-1} \in TM$ and $u\in E$, or $p+k-2=n-1$. The first case contradicts
the minimality of $k$, and the second case contradicts the inequality
$k\le n-p$. This shows that $k=1$, thus proving our assertion. 
\qed

\section{Killing forms on symmetric spaces}

Let $u\in\L^pM$ be a Killing form on a compact simply connected 
symmetric space $M$.

\begin{epr} The pair $(u,du)$ is a section of the bundle $E\oplus F$
defined in the previous section. 
\end{epr}

\begin{proof}
From (\ref{k1}) and (\ref{k2}) we see that $(u,du)$ is a section of
$E_0\oplus F_0$. Moreover (\ref{killing}) and (\ref{killing2}) show
that if $(u,du)$ is a section of $E_k\oplus F_k$ for some $k\ge 0$,
then it is also a section of 
$E_{k+1}\oplus F_{k+1}$. A simple induction argument ends up the proof.
\end{proof}

Suppose now that $p\ge 2$. Then $M$ is not K{\"a}hlerian:

\begin{elem}
A Killing $p$--form on a compact K{\"a}hler manifold is parallel if
$p\ge2$.
\end{elem}

\proof
We make use of the classical K{\"a}hlerian operators 
$$d^c:=\sum
Je_i\wedge\n_{e_i},\qquad \delta^c:=-\sum
Je_i\i\n_{e_i},$$
$$J=\sum Je_i\wedge e_i\i,\qquad \L:=\frac12\sum e_i\i Je_i\i$$
acting on forms, which satisfy the well-known relations 
$$[d,J]=-d^c,\qquad [d,\L]=\delta^c,\qquad d^c\delta+\delta
d^c=0=\delta^c\delta+\delta \delta^c$$ 
on K{\"a}hler manifolds.

Let $u$ be a Killing $p$--form.
We take the wedge product with $JX$ in (\ref{killing}) and sum over an
orthonormal basis $X=e_i$ to obtain:
$$ d^c
u:=Je_i\wedge\n_{e_i}u=\frac{1}{p+1}Jdu=\frac{1}{p+1}(dJu+d^cu),$$
whence
$$p d^c u=dJu.$$
Since $\delta u=0$ and $d^c$ anti--commutes with $\delta$, the two
members of this
relation are $L^2$--orthogonal, so they both vanish.  

Taking now the interior product with $JX$ in (\ref{killing}) and summing over an
orthonormal basis $X=e_i$ yields
$$\delta ^c u=-Je_i\i\n_{e_i}u=\frac{2}{p+1}\L
du=\frac{2}{p+1}(d\L u+\delta ^c u),$$
so 
$$(p-1)\delta ^c u=d\L u.$$
By $L^2$--orthogonality again, (using also the fact that $p\ge 2$) we
get $\delta ^cu=0$. Thus $\Delta u= d^c\delta ^cu+\delta ^c d^cu=0$,
showing that $du=0$, and by (\ref{killing}), $\n u=0$.

\qed

This result shows that we can assume that $M$ is not K{\"a}hlerian. If $M$ is
reducible, then $u$ is a sum of pull--backs of Killing forms on the
factors (see \cite{prod}). 
We can therefore suppose, without lost of generality, 
that $M$ is irreducible. From
Corollary (\ref{cor2}) we deduce that either $R^+(X)u=0$ for every $X$,
or $E=\L^pM$. The first case implies that $du$ is parallel, so in
particular $\Delta u=0$. The Weitzenb{\"o}ck formula (cf. \cite{uwe})
$$\Delta u=(p+1)\n^*\n u$$
then shows that $u$ is parallel.

Consider now the second possibility: $E=\L^pM$. 

\begin{elem} If (\ref{k1}) holds for any $p$--form $u$, $2\le p\le n-2$, 
then the Weyl tensor of $M$ vanishes.
\end{elem}
\proof The equation (\ref{k1}) is $O(n)$--invariant, so if it holds for 
a given non-zero curvature tensor, it must hold for all curvature tensors 
belonging to the corresponding $O(n)$--invariant space.

The equation holds trivially for the scalar part of the curvature tensor, 
and, on the other hand, all symmetric spaces are Einstein, so we only
need to consider Ricci--flat curvature tensors. Therefore, 
if (\ref{k1}) holds 
for the curvature tensor $R$, it equally holds for $W$, where $W$ is the 
Weyl component of $R$. If this is non--zero, (\ref{k1}) must hold for all
tensors of Weyl type, because the space of Weyl tensors is 
$O(n)$--irreducible for $n\ge4$. (For $n=4$ there are two 
$SO(4)$--irreducible 
components, but these are distinguished by the orientation only, so 
they are not $O(4)$--invariant).

We will give an example of a Weyl tensor in dimension 4, and of a
particular 
$2$--form $u_0$, for which  (\ref{k1}) fails.

In higher 
dimensions we complete this Weyl tensor in the trivial way, and for 
higher degree forms we simply take products of $u_0$ with some $p-2$ 
form depending only on the last $n-4$ variables. Note that this operation 
will produce examples of Weyl tensors and $p$--forms that do not 
satisfy (\ref{k1}), as long as $2\le p\le n-2$.

Consider $\a:=g(I\cdot,\cdot), \b:=g(J\cdot,\cdot), \g:=g(K\cdot,\cdot)$
 a basis of self--dual $2$--forms $\L^+$ in $\RM^4$, obtained by composing the
Euclidean metric with three orthogonal complex structures on $\RM^4$
that induce the quaternionic structure on $\RM^4=\HM$. In suitable
 coordinates we have 
$$\ba{ccccc}\a&=&e_1\wedge e_2&+&e_3\wedge e_4,\\
\b&=&e_1\wedge e_3&-&e_2\wedge e_4,\\
\g&=&e_1\wedge e_4&+&e_2\wedge e_3.\ea$$

Define the curvature tensor $R$ by $R(\a)=\a\simeq I,\ R(\b)=-\b\simeq -J, 
R(\g)=0$, 
and extend it by $0$ on anti--self--dual $2$--forms. It
is a Ricci--flat, self--dual curvature tensor.

Let us compute $R^+(X)u$, for any $2$--form $u$:
$$\rp(X)u=e_i\wedge R_{X,e_i}u=e_i\wedge\left(\frac12\langle X\wedge e_i,\a\ra 
I-\frac12\langle X\wedge e_i,\b\ra J\right)u,$$
where the factors $1/2$ come from the fact that $\|\a\|^2=\|\b\|^2=2$.

Setting $u=u_0:=\b$ and using $I\b=2\g$ and $J\b=0$, we get
$$\rp(X)\b=e_i\wedge\g\la X\wedge e_i,\a\ra=e_i\wedge\g\la X,I(e_i)\ra
=-I(X)\wedge \g=J(X)\i\o,$$
where $\o=e_1\wedge e_2\wedge e_3\wedge e_4$ is the volume form.

We get 
$$e_1\i\rp(e_2)\b-e_2\i\rp(e_1)\b=\g,$$
 but $R_{e_1,e_2}\b=0$, which contradicts (\ref{k1}).
\qed

So, if $E=\L^p(M)$, and the locally symmetric space $M$ is not a space form, 
then $p$ must be either $1$ --- in which case we get the Killing vector 
fields --- or $p=n-1$.

In this latter case, the Hodge dual --- say $\xi$ --- of $u$ is a 
closed $1$--form satisfying the equation
\beq\label{fin}\n_X\xi=-\frac1n \delta\xi X.\eeq
In particular $d\xi=0$, and since $M$ is Einstein, the Bochner formula yields 
$$n\xi=\Ric(\xi)=\Delta\xi-\n^*\n\xi=\frac{n}{n-1}d\delta\xi,$$
so $\Delta(\d\xi)=(n-1)\delta\xi$. If $\d\xi=0$, (\ref{fin}) shows
that $\xi$ -- and thus $u$ -- is parallel. Otherwise, the Obata
theorem (cf. \cite{obata}) shows that
$\delta\xi$ is a characteristic function of the round sphere, so
the Weyl curvature of $M$ vanishes. This proves Theorem \ref{main}.


 \labelsep .5cm

\end{document}